\documentclass{llncs}
\usepackage[T1]{fontenc}
\usepackage{xcolor}
\usepackage{bm}
\usepackage{graphicx}
\usepackage{amssymb}

\usepackage{comment}

\usepackage[
    backend=biber,
    style=apa]{biblatex} % <- AGREGAR ESTA LÍNEA
\begin{document}
%
% ENGLISH VERSION
\title{On the Structural Limits of Machine Learning Decision Systems: An Information-Theoretic, Interaction-Based, and Stochastic-Dynamical Perspective}

\author{Nestor Ruben Barraza\inst{1,2}\orcidID{0000-1111-2222-3333} \and
Gabriel Pena\inst{1,2}\orcidID{1111-2222-3333-4444}}
%Third Author\inst{3}\orcidID{2222--3333-4444-5555}}

\authorrunning{N. Barraza et al.}

\institute{Universidad Nacional de Tres de Febrero. Departamento de Ciencia y Tecnología. \and
Universidad de Buenos Aires. Facultad de Ingeniería.}

\maketitle              % typeset the header of the contribution

\begin{abstract}

Machine learning procedures are commonly evaluated in terms of predictive accuracy and computational efficiency. However, their achievable performance is fundamentally constrained by structural properties of the underlying data-generating process, which are formalized in terms of informational bounds. In this work we examine intrinsic limits of data-driven decision systems from an information-theoretic and interaction-based perspective. We analyze minimal achievable error in classification through Fano-type bounds and precision limits in parametric estimation via the Cramér-Rao inequality, emphasizing that such limits depend on the underlying model rather than on algorithmic sophistication alone. We further discuss how implicit assumptions, such as independence, ergodicity, and distributional stability, affect the validity of inferential procedures. Building on interaction-based modeling principles, we review typical frameworks such as Markov Random Fields and potential-based representations for encoding dependence mechanisms. We also describe decision systems, including LLM-integrated agent architectures, as feedback-driven stochastic processes where state-dependent dynamics may induce emergent macroscopic behavior. This perspective highlights the importance of having adequate models for the data as a prerequisite for expanding predictive capability, and situates algorithmic learning within the informational limits imposed by the models.

\keywords{machine learning\and information theory\and complex systems\and stochastic modeling}
\end{abstract}

% Spanish title and abstract
\newpage

\begin{center}
\vspace{1cm}
{\Large \bfseries\boldmath\ Sobre los límites estructurales de los sistemas de decisión mediante aprendizaje automático: una perspectiva desde la teoría de la información, las interacciones y la dinámica estocástica}

\vspace{0.5cm}
\end{center}

\begin{abstract}
Los procedimientos de aprendizaje automatizado usualmente se evalúan en términos de poder predictivo y eficiencia computacional. Sin embargo, la \textit{performance} alcanzable se encuentra limitada por propiedades estructurales del fenómeno subyacente, que se formalizan en términos de cotas de información. En este trabajo se examinan los límites intrínsecos desde una perspectiva de la teoría de la información y las interacciones. Analizamos el mínimo error alcanzable en problemas de clasificación a través de cotas tipo Fano y los límites de precisión en la estimación paramétrica a través de la desigualdad de Cramér-Rao, enfatizando que dichos límites dependen únicamente del modelo de los datos y no del algoritmo utilizado. Adicionalmente, discutimos como suposiciones implícitas tales como independencia, ergodicidad y estabilidad distribucional afectan la validez de los métodos de inferencia. Trabajando sobre principios de modelado basados en interacciones, revisamos \textit{frameworks} clásicos como los campos aleatorios de Markov y representaciones basadas en potenciales para codificar dependencias complejas. También describimos los sistemas de decisión, incluyendo arquitecturas de agentes integrados con modelos de lenguaje, como procesos estocásaticos realimentados donde la dinámica dependiente del estado induce comportamientos emergentes. Esta perspectiva destaca la importancia de disponer de modelos adecuados para los datos como un prerrequisito para mejorar la capacidad predictiva, situando así los algoritmos y métodos de aprendizaje dentro de las limitaciones teóricas fundamentales que los gobiernan.

% palabras clave en minúsculas
\keywords{aprendizaje automatizado\and teoría de la información\and sistemas complejos\and modelado estocástico}
\end{abstract}

\vspace{0.5cm}
\setcounter{page}{1}

\section{Introduction}

In the recent years, Machine Learning (ML) has achieved remarkable success across a wide range of domains such as vision, language processing, finance, logistics, and decision support systems. Advances in optimization, model architectures and large-scale data processing have significantly improved predictive performance in many practical applications. Despite these achievements, predictive accuracy and computational efficiency alone do not suffice to enhance the capabilities of learning systems. Any data-driven method operates under structural constraints imposed by the underlying generative mechanism. Observed data do not constitute an unlimited source of information, and the uncertainty cannot be reduced arbitrarily.

From this perspective, it becomes essential to distinguish between algorithmic sophistication and structural informational limits. Increasing model complexity or data volume does not eliminate intrinsic uncertainty, nor does it alter the informational content encoded in the joint distribution of the involved variables. All ML procedures have structural constraints, which can be formalized by a collection of tools from statistics and information theory. We focus mainly on two cases. For classification tasks, the uncertain is measured through the conditional entropy, and the error probability is bounded by Fano's inequality. For parametric regression methods, the attainable precision is bounded by the Fisher information, through the Cramér-Rao relation.

Moreover, many real-world systems exhibit complex interactions that defy classical models. Reinforcement, long-range dependence, and distributional instability, among others, contradict the usual assumptions and may influence the outcome of procedures that requires them. Thus learning capabilities can be seriously limited by two factors: algorithm independent informational bounds and algorithm specific misassumptions.

The common thread underlying the ideas presented in this work is that predictive performance cannot be fully understood from learning algorithms alone. Information-theoretic limits, interaction models, and stochastic-dynamical descriptions are examined as complementary tools for characterizing the generative mechanisms that produce the data. From this perspective, understanding the phenomenon behind the observations becomes a prerequisite for understanding the capabilities and limitations of any learning procedure.

In this work, we propose a conceptual framework to analyze ML decision systems through two complementary lenses: information-theoretic limits and interaction-based modeling. First, we formalize intrinsic bounds on achievable performance in classification and regression tasks. Following, we emphasize the role of interactions in shaping informational quantities and emergent behavior. Finally, we describe modern decision architectures, including large language models (LLM) agent systems, as stochastic dynamical processes with feedback. This viewpoint connects structural modeling with path dependence and emergent macroscopic behavior, providing a unified perspective on the limits and possibilities of data-driven decision systems.

\section{Structural restrictions}

In the latest years, researchers have focused thoroughly on improving existent algorithms and developing new ones. However, not that much effort has been invested on data modeling. We want to emphasize how adequate modeling is strictly necessary, as model-less algorithms have structural limitations. Information theory provides an important conceptual foundation for understanding data‐driven systems and ML pipelines by relating uncertainty, information content and predictive performance~\parencite{sage}. 

\subsection{Classification problems}

Information cannot be created from nothing. The outcome of any system or algorithm is ultimately constrained by the information coded in its input. Thus, achievable performance in any learning problems is not infinite; there exist informational bounds that are algorithm independent.

Modern learning theory often focuses on algorithmic issues such as optimization strategies and computational performance. This is, however, ussually done without any consideration about the mentioned structural limitations. Algorithms extract and reorganize the information already present in the data-generating process; if there is no such information, then no algorithm can learn anything, no matter how sophisticate. 

There are two fundamental results that formalize these statements. The first is the data processing inequality, which states that no processing can create information, i.e.:
\begin{equation}
    \label{dp_inequality}
    I(X,\hat{X}) \leq I(X,Y),
\end{equation}
where it is assumed that $X \rightarrow Y \rightarrow \hat{X}$ forms a Markov chain. Here $X$ denotes an unknown random variable (r.v.), $Y$ an observed r. v., $\hat{X}$ an estimator for $X$ built on top of $Y$ and $I(\cdot, \cdot)$ the mutual information~\parencite{Scarlett_Cevher_2021}. The inequality implies that any attempt to estimate $X$ through $Y$ cannot create any new information, thus imposing an upper bound in the achievable performance of any algorithm. This bound can be raised if we have more information available in $Y$; for example, if $Y$ denotes a model of the data, a bad model would lead to a tighter upper bound. Thus constructing adequate models is directly related to learning performance.

The second key result is Fano inequality, which is often written in the general form:
\begin{equation}
    \label{fano}
    H(X | \hat{X}) \leq H_2(P_e) + P_e \log(|\mathcal{X}| - 1),
\end{equation}
with the assumption that $X$ is a discrete random variable with support on the set $\mathcal{X}$ (see e.g.~\cite{Cover_Thomas,Scarlett_Cevher_2021}). Here $H(X|\hat{X})$ is the conditional entropy, $H_2$ denotes the binary entropy function. $P_e$ is the error probability, which depends on the problem formulation (i.e., what do we call an error). For a classification problem with exact recovery, we can write $P_e = P(X \neq \hat{X})$, where $X$ and $\hat{X}$ represent the correct and estimated classes. By writing $H(X | \hat{X}) = H(X) - I(X,\hat{X})$ and bounding $H(X)$ by its worst case (i.e., when $X$ is uniform) we have:
\begin{equation}
    \label{fano_2}
    P_e \geq 1 - \frac{I(X, \hat{X}) + \log 2}{\log|\mathcal{X}|}.
\end{equation}
This imposes a bound on the error probability given by the mutual information between the variables. Note that if the support $\mathcal{X}$ is not finite but the mutual information is finite, then we will see an error almost surely. For a finite support (often called ``alphabet'' in this setting), the bound is proportional to the mutual information. Thus, again, a worse choice of model (or not a choice at all) will make any learning system unable to predict better than a fixed upper bound. The approximate recovery situation (where sufficiently close estimated values are accepted) is similar, and is also treated in~\cite{Scarlett_Cevher_2021}. 

\subsection{Regression problems}

In a regression problem, the aim is to provide estimates about some unknown continuous r.v. $X$. The classical formulation of this problem implies finding the parameter vector $\bm{\theta}$ that minimizes some loss function $L$ with respect to some specified parametric family $p(X, \bm{\theta})$. Any learning task that can be formulated in this manner is subject to the Cramér-Rao bound:
\begin{equation}
    Cov_{\bm{\theta}} (\bm{\hat{\theta}}) \succeq { \phi(\bm{\theta}) \, \mathcal{I}({\bm{\theta}})}^{-1}  \phi(\bm{\theta})^T,
\end{equation}
with $\phi(\bm{\theta}) = I + \nabla_{\bm{\theta}}\, b(\bm{\theta})$. The symbol $\succeq$ denotes semidefinite positive order, $\mathcal{I}$ is the Fisher information matrix and $\bm{\hat{\theta}}$ is an estimator with bias function $b(\bm{\theta})$. If the estimator is unbiased then the lower bound reduces to the inverse of the Fisher information matrix.

As in the classification case, the learning performance (described here in terms of the covariance matrix of the estimator) is constrained by the information encoded in the data. This is not algorithmic but structural: estimation cannot be arbitrarily precise. 

The most important conclusion to obtain from this section is the following: \textbf{informational bounds are structural and algorithm independent}. The only way to overcome them is to choose suitable models to better describe the data.

\section{Statistical limitations}

In the previous section we discussed informational constraints that were algorithm independent. We now turn to algorithmic-specific constraints, which are often underestimated or not considered at all. These arise mostly due to implicit statistical assumptions that do not hold.

\subsection{The Central Limit Theorem}

Most supervised learning frameworks make implicit assumptions regarding correlations between observations. The simplest case is when observations are assumed to be a random sample, i.e., independent and identically distributed (IID). Under this assumption one has guaranteed almost-sure (a.s.) convergence from the Law of Large Numbers (LLN) and convergence in distribution from the Central Limit Theorem (CLT), which justify most of the classical inference procedures. The CLT can be generalized to allow certain forms of dependence among observations; in broad terms, such results require finite second moments and a sufficiently weak dependence structure. Typical conditions include summable autocorrelations or mixing assumptions, which ensure that correlations decay sufficiently fast as the lag increases. These conditions are not always guaranteed. If any of those fail, the resulting dynamics may show phenomena known in physics as \textit{anomalous diffusion}, where the scaling differs significantly from the expected Gaussian behavior.

Infinite variance is a staple property of heavy-tailed distributions. In the most extreme cases, even the first moment may diverge; a well-known example is given by Pareto-type distributions. These distributions are problematic in the sense that they fail to have the expected statistical behavior. If the mean value is infinite then averages are dominated by rare events; if the mean is finite but the variance is not, then large observations do not dominate but occur occasionally. As a consequence, the usual Gaussian approximation may not be adequate to describe the limit distribution. On the other hand, joint distributions where correlations decay slowly are also troublesome, since results that assume approximate independence cannot be applied.  

A special case arises when the observations come from a stochastic process $(X_n)_n$ indexed by $n$. If the samples are interpreted as increments or displacements and one wants to address the distribution of the total motion distance $X^{(n)} = \sum_{k=1}^n X_k$, then for the CLT to hold an additional condition has to be met: strict stationarity~\parencite{VilkAnomalousDiffusion2022}, i.e., the joint distribution being invariant with respect to index shifts. Thus, aggregated processes $X^{(n)}$ with increments that are not strictly stationary are not guaranteed to converge to a Brownian motion. 

Also, in the case of stochastic processes, slowly decaying correlations induce memory effects that persist over long time scales; this phenomenon is called long-range dependence (LRD). A classical example is fractional Brownian motion with Hurst exponent $H\geq1/2$, where the correlation decays as a power law. In this case, LRD prevents early fluctuations from dissipating rapidly, and asymptotic behavior can differ significantly from the Brownian case. This highlights that knowledge of the underlying correlation structure should not be neglected, as classical results do not always apply. It is also worth noting that LRD does not contradict Markovianity; this is a property concerning conditional transition probabilities, whereas LRD is a statement about moments. As an example, the 3p-BPM$(\beta, \gamma, \rho)$ process~\parencite{Barraza_2025,CSFEpidemiologia} has LRD whenever its parameters satisfy $\gamma/\rho \neq 1/2$, yet remains Markovian~\parencite{fendick_whitt_2022,Barraza_2025}.

\subsection{Ergodicity}

Ergodicity means, roughly speaking, that time averages approximate ensemble averages. In other words, having a single realization of a stochastic process split into short segments and averaged is equivalent to having several short trajectories. Although many well-known processes are ergodic (e.g. $AR(1)$ autoregressive models, the Ornstein-Uhlenbeck process and the Poisson process, among others) it is not nearly universal, and most processes found in nature are not ergodic. 

In practice, it is often impossible or expensive to attain several sample trajectories of a process, but it is acceptable to observe a single realization for a long period. If the process is ergodic then one can safely perform statistical inference by averaging observations in time; but this need not be the case. Blindly applying methods that require ergodicity will thus not provide reliable results. 

Non-ergodicity is closely related to the persistence of structural components across time: different realizations may converge to different time averages. This means that asymptotic statistics can depend on the specific invariant subspace selected by the initial state of the system, i.e., that the influence of the trajectory starting point may never dissipate. From the algorithmic point of view, this implies that any learning or inference method that relies on a single long trajectory may not reflect ensemble statistics, since the observed realization may belong to a specific statistical regime that need not be the only possible one.

\subsection{Frequency estimators}

Learning models often seek to optimize a function $L(\bm{x}, \bm{y}, \bm{p})$, where $\bm{x}$ denotes an input vector and $\bm{y}$ denotes an output vector, with respect to a parameter vector $\bm{p}$, which may have very high dimensionality (as in the case of deep neural networks). The usual method involves minimizing a loss function, for which there exist many options, such as cross-entropy optimization, mean square error ($MSE$) and penalized methods such as Ridge or Lasso.  

In explicit parametric models such as neural networks, the optimization procedure is performed globally, in the sense that one seeks to minimize the divergence between softmax outputs and empirical frequencies by exploring the parameter space (i.e., through gradient descent or annealing methods). There is no explicit estimation of the probabilities. In contrast, many models (mostly non-parametric) require such an explicit estimation. In these settings, probabilities are estimated by counting frequencies over pre-defined regions of the feature space:
\begin{equation}
    \hat{P}(\bm{Y} = \bm{y} | \bm{X} = \bm{x}) = \frac{\# \mbox{ positives}}{\# \mbox{ total count}}.
\end{equation}
Examples include naive Bayes, k-Nearest Neighbors, histograms and classical $n-$grams language models. Explicit estimation of probabilities via frequencies has some limitations. Convergence of the empirical frequency to the real probability is guaranteed by the LLN, provided its hypothesis hold ---which as we have discussed, does not always happen---. Therefore a first restriction is given by the dependence structure of the data. In models with high dimensionality or continuous range point probabilities are always zero, so one defines regions in the feature space. The typical example is histograms. Even though there are formulas and methods to define regions according to some criteria, this separation is mostly \textit{ad-hoc}, which induces biases and variability. Another classical example where this estimation fails is $n-$grams, where the samples are usually sparse (i.e., most scenarios never or almost never happen), leading to systematic underestimation of some events and overestimation of others. This problem can be overcome by using corrections or alternative estimators such as Good-Turing, Kneser-Ney or Laws of Succesion, among others~\parencite{10.1063/1.3038987}.

\section{Modeling}

A recurrent theme in contemporary ML is the search for increasingly expressive hypothesis classes and better optimization procedures. Yet algorithmic development, as we have discussed, is not sufficient. We shall address a more fundamental question: how are interactions between variables represented?

\subsection{Interaction-based modeling}

Any learning model encodes, explicitly or not, assumptions about interactions. For example, the linear regression model assumes that variables interact additively (mediated by coefficients), whereas neural networks use activation functions to introduce non-linear interactions. In a general setting, modeling can be thought as specifying how different components of a system influence each other. These influences are called \textit{interactions}, and the corresponding modeling techniques are called \textit{interaction-based modeling}. 

In many real world systems, independence assumptions are structurally inadequate, since observations can be spatially or temporarily correlated, socially reinforced or coupled, among several possible interactions. In the probabilistic framework, interactions are encoded in the joint distribution of the features. Informational quantities such as the Fisher information and the entropy derive directly from this. Therefore, adequately choosing probabilistic models to represent these dependencies is a fundamental problem to be addressed.

Markov Random Fields (MRFs) provide a formal framework to represent interactions with an intrinsic geometric structure. Let the random vector $\bm{X} = (X_1, \dots, X_d)$ describe the complete state of a system, where each index corresponds to a vertex in an undirected graph $G=(V, E)$, where each edge in $E$ encodes an interaction. A MRF is a strictly positive distribution satisfying the local Markov property with respect to $G$ (i.e., vertexes are conditionally independent given its neighbors). The Hammersley-Clifford theorem establishes that the joint distribution of $\bm{X}$ factorizes over the cliques of $G$~\parencite{hammersley1971markov,besag1974spatial,griffeath1976mrf}:
\begin{equation}
    P_{\bm{X}}(\bm{x}) = \frac{1}{Z} \exp\left(- \sum_{C \in \mathcal{C}} \phi_C(x_C)\right),
\end{equation}
where $\phi_C$ are clique potentials, $\mathcal{C}$ is the set of cliques $C$, $x_C=\lbrace x_i; \, i\in C \rbrace$ is the state vector restricted to $C$, and $Z$ is the normalization constant
\begin{equation}
    Z = \sum \exp\left(- \sum_{C\in\mathcal{C}} \phi_C(x_C)\right)
\end{equation}
which in statistical mechanics is called \textit{partition function}.

Clique potentials represent local interactions, which combined and subject to a proper normalization describe global properties. Self-potentials to capture reflexive relations in a single vertex are also acceptable. In this case, the modeling task reduces to capturing structural properties by meaningful potentials. Table \ref{tab_mrfs} briefly describes some examples. 

\begin{table}[!htb]
    \centering
    \begin{tabular}{l l p{6cm}}
    \hline
    Interaction & Potential $\phi_{ij}(x_i,x_j)$ & Interpretation \\
    \hline
    Alignment
    & $-\beta x_i x_j \qquad \beta > 0$ 
    & Encourages similarity or agreement. Captures positive coupling such as contagion phenomena (e.g. infectious diseases or herding in markets, Ising model in statistical mechanics).\\
    \hline
    Anti-alignment 
    & $+\beta x_i x_j \qquad \beta > 0$ 
    & Promotes differentiation or competition by penalizing similarity. Can model frontiers or boundaries in image segmentation or product substitution in market analysis.\\
    \hline
    Spatial smoothness 
    & $\lambda (x_i-x_j)^2 \qquad \lambda > 0$ 
    & Enforces local continuity, promoting smooth transitions among regions. Induces a Gaussian MRF (see e.g. \cite{Rasmussen_Williams_Gaussian}, Appendix B5). Usage examples include environmental measurements, spatial variations in demand and audio/image smoothing. \\
    
    \hline\hline
    
    Interaction & Self-potential $\phi_i(x_i)$ & Interpretation \\
    \hline
    Reinforcement 
    & $\propto -\alpha \log(x_i + c) \qquad \alpha > 0$ 
    & Accumulative inertia. Marginal feature distribution exhibits a heavy-tailed law, thus admitting extreme outliers with non-negligible probability. Captures individual accumulation dynamics, such as wealth or popularity. Resembles Polya-type processes. \\
    \hline
    Capacity
    & $\gamma x_i^2 \qquad \gamma>0$ 
    & Introduces limitations in capacity by producing Gaussian marginal distributions centered around zero. Can be used to model saturation and congestion in supply-chain or traffic systems. \\
    \hline
    
    \end{tabular}
    \caption{Some examples on modeling interactions via clique potentials.}\label{tab_mrfs}
\end{table}

The partition function integrates local interactions into a global joint distribution structure, which in turn determines systemic properties. A classical example is the Ising model~\parencite{MezardInformationPhysics,BaxterStatisticalMechanics}, where spontaneous magnetization appears as an emergent effect from pairwise spin interactions. Similar models have been successful in explaining large-scale organization and observed macroscopic states in statistical mechanics. In data-driven contexts, interaction based models such as MRFs provide a structured way to encode domain knowledge in the joint distribution of system states.

Simple local interaction rules can generate complex global phenomena. This is a remarkable consequence of interaction-based modeling, observed in both natural and artificial contexts. The Ising model is such an example. Neural networks operate under the same principle: each neuron encodes very little information, through the set of weighted connections, but an enormous number of neurons can produce impressive results, be it creating approximate solutions for the traveling salesman problem via a self-organizing map~\parencite{Kohonen2001} or achieving powerful generative capabilities~\parencite{Vaswani2017AttentionIA}. In all these examples, individual behavior or local interactions determine collective properties.

\subsection{Learning and decision in path-dependent multi-agent systems}

Modern data-driven systems ---including multi-agent models, usually coupled with one or more large-language models (LLMs)--- can be seen as stochastic dynamical systems that evolve through time, subject to feedback and internal interactions. Rather than acting as a static deterministic mapping, these systems keep internal states that accumulate information to influence the future behavior, which is inherently stochastic~\parencite{4445757,shukla2025stochasticdifferentialequationframework,carson2025stochasticdynamicaltheoryllm}. In particular, agent-based systems governed by macroscopic physical laws have been recently proposed in~\cite{song2025detailedbalancelargelanguage}. In the same work, the existence of ``microscopic'' statistical interaction laws is observed by measuring detailed balance on an experimental framework.   

A general formulation for these systems can be made by the triple $X_t = (S_t, A_t, O_t)$, where the index $t$ denotes time (either continuous or discrete). The variable $S_t$ represents the current system state (encoding accumulated information), $A_t$ represents an action following some decision rule and $O_t$ represents an observation from the environment. The dynamics evolve as follows: the next state $S_{t+1}$ is determined as a function of the current configuration $X_t$, actions are sampled from a decision policy (which may be random and depend on the current state), and observations are drawn from a distribution that depends on previous samples (environment) and the chosen action. 

In this framework, internal feedback arises when states or actions depend on the current state $S_t$, which thus encodes historical information from previous interactions. This mechanism induces path-dependence, in the sense that the future depends (directly or indirectly) on the process internal history. In a Markovian setting, the dependence occurs entirely through the state variable $S_t$. We remark that Markovianity does not contradict path-dependence: it just imposes a specific mode in how this dependence is structured (i.e., conditionally independent of the past given the present).

In a multi-agent setting, the state $S_t$ can be interpreted as a system-wide configuration of several interacting agents and environmental inputs. Each agent produces actions according to its own decision rule, which may be conditioned on local observations or internal memory; the global behavior is thus an emergent property that arises from the interaction between those agents. From a modeling perspective, this system can still be represented as a single stochastic dynamical process by letting the state space take into account all the internal variables, from all agents (see e.g. \cite{4445757}, where a multi-agent system is modeled as a stochastic game). This formulation  provides a mathematically tractable framework to describe complex systems in terms of one or more interaction models. It also makes explicit that path-dependence arises naturally from the way agents interact with each other (e.g. providing feedback).

\subsection{Path-dependent Markovian dynamics}

In classical macroscopic continuous descriptions, dynamical processes $X_t$ are described by a stochastic differential equation (SDE) of the form:
\begin{equation}
    dX_t = \mu(X_t, t) \, dt + \sigma(X_t, t) \, dW_t,
\end{equation}
where $\mu$ and $\sigma$ are parameter functions and $W_t$ denotes a Wiener process (Brownian motion). The induced probability density then evolves according to a Fokker-Planck equation:
\begin{equation}
    \frac{\partial p(x, t)}{\partial t} = -\frac{\partial}{\partial x} [\mu(x, t) p(x, t)] + \frac{1}{2} \frac{\partial^2}{\partial x^2} [\sigma(x, t) p(x, t)].
\end{equation}
This equation is often interpreted in the following way: the process evolves through time with a drift given by $\mu(x, t)$ and a diffusion behavior determined by $\sigma(x, t)$. In particular, if the drift and diffusion coefficients do not depend on the position $x$ then the process does not have a structural feedback dynamics. If any of the parameter functions depend on $x$, then it encodes a feedback mechanism and the evolution becomes state-dependent. Since all the relevant historical information is embedded into the present state $X_t$, this class of processes remain Markovian, despite being possibly path-dependent.

While this formulation produces a useful macroscopic description of continuous processes, many natural and man-made systems can be better explained through discrete state-space models. These cannot be directly obtained by the Fokker-Planck formulation; instead, they require to be formulated in terms of counting processes. A wide class of such are the \textit{birth and death processes}, governed by the following infinite system of differential equations:
\begin{eqnarray}
    &\frac{d}{dt} p(k, t) = -[\lambda(k, t) + \mu(k, t)]\, p(k, t) + \nonumber\\
    &+ \lambda(k-1, t)\, p(k-1, t) + \mu(k+1, t)\, p(k+1, t) \qquad k\geq 1,\\
    &\frac{d}{dt} p(0, t) = -[\lambda(0, t) + \mu(0, t)] p(0, t) + \mu(1, t) p(1, t)
\end{eqnarray}
where $k \in\mathbb{N}$ and $t > 0$ and $\lambda(k,t), \mu(k,t)$ are functions called birth rate and death rate respectively. These are the discrete state-space analogous of the Fokker-Planck: both correspond to the Kolmogorov forward equation, in different contexts. These equations govern the probability of a particle performing a transition of size $k$ in a time interval of length $t$. Feedback is encoded in the same manner: if either the birth or death rate depends explicitly on the system state $k$, then the dynamics become path-dependent. A typical example is given by the class of processes with rates of the form
\begin{equation}
    \lambda(k, t) = (\lambda_0 + \alpha k) \kappa(t) \qquad\qquad \mu(k, t) \equiv 0,
\end{equation}
called Generalized Polya Processes (GPP). The GPP describe situations in which the only allowed transitions are jumps of unit length (pure birth processes). In a GPP, the birth rate is proportional to the position $k$, thus encoding a positive reinforcement or contagion~\parencite{Barraza_2025,PenaMorenoBarraza2022,CSFEpidemiologia}. GPP can be also seen as the continuous-time counterparts of the Polya urn model and its generalizations~\parencite{urn_models}. This illustrates how structural modeling of the interactions reshapes the global behavior of the system, generating non-trivial emergent properties.

\section{Conclusion}

In this work, we have examined how machine learning systems operate within structural constraints determined by the underlying generative phenomenon. We showed that informational quantities such as conditional entropy and Fisher information establish fundamental bounds on achievable classification accuracy and estimation precision. These limits are intrinsic to the data-generating process and algorithm independent, since they depend only in the information encoded within the data. We emphasized that implicit assumptions ---including independence, correlation structure and ergodicity, among others--- play a decisive role in the effectiveness of algorithms. In consequence, performance degradation can arise from a structural specification mistake rather than from purely algorithmic or computational limitations.

The importance of understanding the data-generative phenomenon rather than merely increasing algorithmic complexity was highlighted throughout this work. In particular, we remarked that improving the predictive capability often requires refining model assumptions and enriching the hypothesis space that encodes the generative structure. Typical interaction-based models such as Markov Random Fields, and dynamical frameworks such as multi-agent systems and diffusive processes were reviewed as useful tools for describing emergent properties ---such as path-dependence--- in terms of structural interactions (reinforcement, feedback), thus providing a relevant framework for machine learning applications.

A central message of this work is that the study of learning systems should not be restricted to the algorithmic layer. Information-theoretic limits, interaction structures, and stochastic dynamics provide complementary perspectives for understanding the mechanisms that generate the data. In this sense, looking beyond the algorithm and towards the underlying phenomenon may be as important as improving predictive architectures themselves.

Viewing decision systems as stochastic dynamical processes clarifies that learning algorithms reorganize information already embedded in the phenomenon. Recognizing structural limits should not be interpreted as a critique of machine learning methodologies, but as a necessary step toward their robust, interpretable, and responsible deployment in complex environments.

\begin{credits}
\subsubsection{\ackname} We would like to thank Universidad Nacional de Tres de Febrero for its support under grant 80020240600005TF.

\end{credits}

\printbibliography

\end{document}